\documentclass[draft]{amsart}
\usepackage{fullpage}

\usepackage{amssymb, amsmath}

\usepackage[english,francais]{babel}

\newtheorem{theorem}{Theorem}[section]

\newtheorem{e-proposition}[theorem]{Proposition}

\newtheorem{e-definition}[theorem]{Definition\rm}

\def\og{\leavevmode\raise.3ex\hbox{$\scriptscriptstyle\langle\!\langle$~}}
\def\fg{\leavevmode\raise.3ex\hbox{~$\!\scriptscriptstyle\,\rangle\!\rangle$}}

\begin{document}
\centerline{Comptes Rendus Math\'ematique. Acad\'emie des Sciences. Paris}


\selectlanguage{english}
\title{Solutions of the Vlasov-Maxwell-Boltzmann system with long-range interactions}


\selectlanguage{english}

\author{Diogo Ars\'enio}
\email{arsenio@math.univ-paris-diderot.fr}
\address{Institut de math\'ematiques, Universit\'e Paris Diderot, Paris}

\author{Laure Saint-Raymond}
\email{saintray@dma.ens.fr}
\address{D\'epartement de math\'ematiques et applications, Ecole normale sup\'erieure, Paris}


\maketitle

\medskip

\begin{center}
{\small Received *****; accepted after revision +++++\\
Presented by £££££}
\end{center}

\begin{abstract}
\selectlanguage{english}
We establish the existence of renormalized solutions of the Vlasov-Maxwell-Boltzmann system with a defect measure in the presence of long-range interactions. We also present a control of the defect measure by the entropy dissipation only, which turns out to be crucial in the study of  hydrodynamic limits.
{\it To cite this article: Diogo and Laure, C. R. Acad. Sci. Paris (2013).}

\vskip 0.5\baselineskip

\selectlanguage{francais}
\noindent{\bf R\'esum\'e} \vskip 0.5\baselineskip \noindent
{\bf Solutions du syst\`eme de Vlasov-Maxwell-Boltzmann avec interactions \`a longue port\'ee. }
Nous \'etablissons l'existence de solutions renormalis\'ees du syst\`eme de Vlasov-Maxwell-Boltzmann avec mesure de d\'efaut en pr\'esence d'interactions \`a longue port\'ee. Nous pr\'esentons \'egalement un contr\^ole de la mesure de d\'efaut par la dissipation d'entropie uniquement, qui s'av\`ere \^etre crucial dans l'\'etude des limites hydrodynamiques.
{\it Pour citer cet article~: Diogo et Laure, C. R. Acad. Sci.
Paris (2013).}

\end{abstract}


\selectlanguage{english}

\section{Introduction}

We consider the Vlasov-Maxwell-Boltzmann system:
\begin{equation*}
	\partial_t f + v \cdot \nabla_x f + \left( E + v \wedge B \right) \cdot \nabla_v f = Q(f,f),
\end{equation*}
where $f(t,x,v)\geq 0$, with $t\in [0,\infty)$, $x\in \mathbb{R}^3$ and $v\in\mathbb{R}^3$, is the particle number density and $(E(t,x),B(t,x))$ is the self-induced electromagnetic field, whose evolution is governed by Maxwell's equations
\begin{equation*}
	\begin{cases}
		\begin{aligned}
			\partial_t E - \nabla_x\wedge B &= - \int_{\mathbb{R}^3} fv dv,
			& \nabla_x\cdot E &=\int_{\mathbb{R}^3} fdv -1, \\
			\partial_t B + \nabla_x\wedge E & = 0, 
			& \nabla_x\cdot B &=0.
		\end{aligned}
	\end{cases}
\end{equation*}
We refer to \cite{arsenio6} and the references therein for the definition of the Boltzmann collision operator $Q(f,f)$ and for a thorough introduction to the mathematical theory of the Vlasov-Maxwell-Boltzmann system. Here, we merely emphasize that $Q(f,f)$ depends on interparticle interactions through a given cross-section $b(z,\sigma)=b\left(|z|,\frac{z}{|z|}\cdot\sigma\right)\geq 0$, where $(z,\sigma)\in\mathbb{R}^3\times\mathbb{S}^2$ and where we are using the $\sigma$-representation of the operator (according to the terminology from \cite{villani}). Assuming that the gas is at equilibrium at infinity (say with density $M(v) = (2\pi)^{-3/2} \exp (-|v|^2/2)$ without loss of generality), the natural a priori bounds on the unknowns $(f,E,B)$ are provided by the relative entropy inequality, for all $t>0$,
\begin{equation}\label{entropy}
	\begin{aligned}
		\int_{\mathbb{R}^3\times\mathbb{R}^3} \left(f \log {f \over M} - f +M\right) (t) dxdv
		& + \frac 1{2} \int_{\mathbb{R}^3} \left( |E|^2+ |B|^2 \right)(t) dx
		+ \int_0^t\int_{\mathbb{R}^3}D(f)(s) dx ds
		\\
		& \leq \int_{\mathbb{R}^3\times\mathbb{R}^3} \left(f^\mathrm{in}\log {f^\mathrm{in} \over M}
		- f^\mathrm{in}+M\right) dxdv
		+ \frac1{2}\int_{\mathbb{R}^3} \left( |E^\mathrm{in}|^2+  |B^\mathrm{in}|^2 \right) dx,
	\end{aligned}
\end{equation}
where the entropy dissipation $D(f)(s)\geq 0$ is formally defined by $-\int_{\mathbb{R}^3}Q(f,f)\log f(s)dv$. Note that the control of the relative entropy implies a bound on $f$ in $L^\infty\left(dt;L^1_{\mathrm{loc}}\left(dx;L^1\left((1+|v|^2)dv\right)\right)\right)$ (see \cite{SR}).

The construction of suitable global solutions to the Vlasov-Maxwell-Boltzmann system for large initial data in the entropy space is considered of outstanding difficulty, due to the lack of dissipative phenomena in Maxwell's equations, which are hyperbolic. Indeed, the DiPerna-Lions construction of renormalized solutions \cite{diperna}, \cite{diperna3}, \cite{lions3} seems to break down as soon as the DiPerna-Lions theory on linear transport equations \cite{diperna6} fails. Thus, so far, the only known answer to this problem is due to Lions \cite{lions3} in the cutoff case, i.e. $b\in L^1_\mathrm{loc}\left(\mathbb{R}^3\times\mathbb{S}^2\right)$, where a rather weak notion of solutions was derived: the so-called measure-valued renormalized solutions. However, these solutions failed to reach mathematical consensus on their usefulness due to their very weak aspect. We wish emphasize that a refinement of such solutions is established in \cite{arsenio6} and that their hydrodynamic limit is shown therein to match standard Leray solutions for the corresponding fluid dynamical systems when they exist. It should be mentioned that an alternative approach yielding strong solutions, provided smallness and regularity assumptions on the initial data are satisfied, was obtained more recently by Guo in \cite{guo}.

In this brief note, we explain how it is possible to establish the existence of renormalized solutions with a defect measure of the Vlasov-Maxwell-Boltzmann system in the presence of long-range interactions, i.e. $b\notin L^1_\mathrm{loc}\left(\mathbb{R}^3\times\mathbb{S}^2\right)$, following the method of Alexandre and Villani \cite{alexandre}. Essentially, we are merely going to explain how the breakdown of the DiPerna-Lions theory on linear transport \cite{diperna6} can be compensated by the regularizing (or at least compactifying) effect due to long-range interactions. This approach provides thus the first global existence result for large initial data in the entropy space (defined by the natural bounds from \eqref{entropy}) for the Vlasov-Maxwell-Boltzmann system. Finally, we note that this method can also be applied to Vlasov-Maxwell-Boltzmann systems of several species of particles and to Vlasov-Maxwell-Landau systems.

\section{Renormalized solutions with a defect measure}

We say that a nonlinearity $\beta\in C^2\left([0,\infty);\mathbb{R}\right)$ is an \emph{admissible renormalization} if it is bounded and satisfies $\left|\beta'(z)\right|\leq \frac{C}{1+z}$ and $\beta''(z)\leq0$, for all $z\geq 0$ and some $C>0$.

A density function $f\geq 0$ and electromagnetic vector fields $E$ and $B$, such that the entropy inequality \eqref{entropy} holds for all times $t>0$, are a \emph{renormalized solution of the Vlasov-Maxwell-Boltzmann system with a defect measure} if $(f,E,B)$ solves for any admissible renormalization
\begin{equation}\label{VMB with defect}
		\partial_t {\beta\left(f\right)} + v \cdot \nabla_x {\beta\left(f\right)}
		+ \left( E + v \wedge B \right) \cdot \nabla_v {\beta\left(f\right)}
		\geq {\beta'\left(f\right)}Q(f,f)
\end{equation}
and Maxwell's equations in the sense of distributions,  and $f$ verifies the continuity equation
\begin{equation*}
	\partial_t\int_{\mathbb{R}^3}f(t,x,v)dv+\nabla_x\cdot\int_{\mathbb{R}^3}f(t,x,v)vdv=0.
\end{equation*}
Furthermore, in order to have a sound notion of weak solution to an initial value problem, we ask that $(f,E,B)$ enjoys some temporal continuity in the weak topology of some appropriate functional space. Note that this definition actually determines a solution and not a mere subsolution (see \cite{alexandre}).

More precisely, the inequality in \eqref{VMB with defect} asks that the difference of the left-hand side with the right-hand side be a positive distribution. Therefore, for any given admissible renormalization $\beta(z)$ there is a positive Radon measure $\nu_\beta$ on $[0,\infty)\times\mathbb{R}^3\times\mathbb{R}^3$ equal to that difference such that $(f,E,B)$ solves
\begin{equation*}
		\partial_t {\beta\left(f\right)} + v \cdot \nabla_x {\beta\left(f\right)}
		+ \left( E + v \wedge B \right) \cdot \nabla_v {\beta\left(f\right)}
		= {\beta'\left(f\right)}Q(f,f) + \nu_\beta
\end{equation*}
in the sense of distributions, thus justifying the appellation of these solutions.

Under the non-cutoff assumption $b\notin L^1_\mathrm{loc}$ (corresponding to long-range interactions in the plasma), it is a priori not clear at all that the right-hand side of the renormalized Vlasov-Boltzmann equation with defect \eqref{VMB with defect} is even well-defined as a distribution. In fact, it has been shown in \cite{alexandre} that the renormalized Boltzmann operator $\beta'\left(f\right)Q(f,f)$ can be defined in some negative index Sobolev space as a singular integral operator, under some very large assumptions on the singular collision kernel $b$ which include all the physically relevant cross-sections, solely based on the bounds provided by the entropy inequality \eqref{entropy}. Other sharper bounds based on the entropy dissipation are given in \cite{arsenio3} and further refined in \cite{arsenio6}, which, in particular, allow for a larger set of admissible renormalizations.

\section{Weak stability of solutions}

Let us consider a sequence of renormalized solutions with defect measure $\left\{\left(f_k,E_k,B_k\right)\right\}_{k\in\mathbb{N}}$ with initial data satisfying the uniform bound
\begin{equation*}
	\sup_{k\in\mathbb{N}}\int_{\mathbb{R}^3\times\mathbb{R}^3} \left(f_k^\mathrm{in}\log {f_k^\mathrm{in} \over M}
	- f_k^\mathrm{in}+M\right) dxdv
	+ \frac1{2}\int_{\mathbb{R}^3} \left( |E^\mathrm{in}_k|^2+  |B^\mathrm{in}_k|^2 \right) dx<\infty,
\end{equation*}
so that the entropy inequality \eqref{entropy} is satisfied uniformly in $k\in\mathbb{N}$. In particular, it is standard procedure to show, up to extraction of a subsequence, using the methods from \cite{diperna3} and assuming the strong convergence of initial data, that $(f_k,E_k,B_k) \rightharpoonup (f,E,B)$, as $k\rightarrow\infty$, at least in $\mathcal{D}'$, where $(f,E,B)$ satisfies the entropy inequality \eqref{entropy}, as well. This first step is a mere lower semi-continuity result on the convex functionals defining the entropy inequality \eqref{entropy} and, thus, is valid under both the cutoff and non-cutoff assumptions from \cite{diperna} and \cite{alexandre}, respectively.

In the non-cutoff setting, which includes all the inverse power kernels (see \cite{villani}), Alexandre and Villani \cite{alexandre} succeeded in showing the strong convergence of the $f_k$'s towards $f$. Their analysis consisted first in exploiting the intricate regularizing (or merely compactifying) effects on the velocity variable due to long-range interactions, thanks to a former study by Alexandre, Desvillettes, Villani and Wennberg \cite{alexandre2} on the entropy dissipation bound. Alternately, one could also use the results from \cite{arsenio5} to deduce a similar velocity smoothness in a more convenient functional setting. Then, they used standard velocity averaging lemmas (see \cite{golse0}, for instance) to obtain some compactness in time and space. This step can be somehow simplified by using the hypoelliptic transfer of compactness obtained in \cite{arsenio}. This approach remains entirely valid for the Vlasov-Maxwell-Boltzmann system and so, we are able to claim here, up to extraction, that $f_k\rightarrow f$ almost everywhere.

The final step of the analysis from \cite{alexandre} established that
\begin{equation}\label{limit defect}
	\liminf_{k\rightarrow\infty}\int_{[0,\infty)\times\mathbb{R}^3\times\mathbb{R}^3}
	\beta'(f_k)Q(f_k,f_k)\varphi(t,x,v)dtdxdv
	\geq
	\int_{[0,\infty)\times\mathbb{R}^3\times\mathbb{R}^3}
	\beta'(f)Q(f,f)\varphi(t,x,v)dtdxdv,
\end{equation}
for every non-negative $\varphi\in C_0^\infty\left([0,\infty)\times\mathbb{R}^3\times\mathbb{R}^3\right)$. In fact, in view of the strong convergence of the $f_k$'s and of the crucial renormalized formulation of the collision operator, this lower semi-continuity is a mere consequence of Fatou's lemma and, thus, still holds for the Vlasov-Maxwell-Boltzmann system.

It is then easy, thanks to the strong compactness of the $f_k$'s, to pass to the limit in the renormalized equation \eqref{VMB with defect} with a defect and to show its stability. Finally, the stability of Maxwell's equations and of the local conservation of mass is easily obtained, for these equations are linear.

It only remains to establish the time continuity of $(f,E,B)$. Even though the temporal continuity of the electromagnetic fields $(E,B)$ in $C\left([0,\infty);\textit{w-}L^2\left(dx\right)\right)$ is readily deduced from Maxwell's equations, it is not as easily exhibited for the density $f$. Indeed, the defect measure $\nu_\beta$ may be singular, even in time, and so, loosely speaking, the renormalized equation \eqref{VMB with defect} only permits us to deduce that $f$ is of bounded variation in time (weakly in space and velocity) and not absolutely continuous. This point is not addressed in \cite{alexandre}. It is discussed in \cite{arsenio3} for the Boltzmann equation on the torus where it is shown that $f\in C\left([0,\infty);\textit{w-}L^1_\mathrm{loc}\left(dxdv\right)\right)$ by exploiting the conservation of mass. A similar argument holds here and we refer to \cite{arsenio6} for a complete justification of time continuity.

On the whole, we have shown that $(f,E,B)$ is a renormalized solutions of the Vlasov-Maxwell-Boltzmann system with a defect measure. In particular, note that, since the strong convergence of the $f_k$'s obviously implied the convergence of the $\beta(f_k)$'s towards $\beta(f)$, there has been no need to appeal to the DiPerna-Lions theory on transport equations \cite{diperna6} as required in the cutoff case (see \cite{lions3}).

The following existence theorem follows from the weak stability of solutions by constructing suitable approximations of the Vlasov-Maxwell-Boltzmann system. Alexandre and Villani elegantly used in \cite{alexandre} actual renormalized solutions to the Boltzmann equation with cutoff to approximate the Boltzmann equation without cutoff. In the present case of the Vlasov-Maxwell-Boltzmann system without cutoff, it is not possible to approximate it with its cutoff analog because standard renormalized solutions are simply not known to exist in this case. Rather, we should use standard truncation and approximation techniques, as in \cite{diperna2, diperna, lions3}, to approximate the Vlasov-Maxwell-Boltzmann system.

\begin{theorem}
	Let $b(z,\sigma)$ be a collision kernel satisfying the non-cutoff assumptions from \cite{alexandre}. Then, for any initial condition $f^\mathrm{in}\geq 0$ and $E^\mathrm{in},B^\mathrm{in}\in L^2\left(dx\right)$, satisfying the initial compatibility conditions $\nabla_x\cdot E^\mathrm{in} =\int_{\mathbb{R}^3}f^\mathrm{in}dv-1$, $\nabla_x\cdot B^\mathrm{in} =0$ and the initial bound
	\begin{equation*}
		\int_{\mathbb{R}^3\times\mathbb{R}^3}\left( f^\mathrm{in} \log {f^\mathrm{in}\over M} -f^\mathrm{in}+M \right)
		dxdv < \infty,
	\end{equation*}
	there exists a renormalized solution $\left(f,E,B\right)$ to the Vlasov-Maxwell-Boltzmann system with a defect measure.
\end{theorem}

We refer to \cite{arsenio6} for a complete justification of the above theorem and for a precise statement of the assumptions on the collision kernel $b$, which are, in fact, slightly weaker than the assumptions from \cite{alexandre}.

\section{A refined entropy inequality}

It is possible to obtain a refined characterization of the defect measure in terms of the entropy dissipation. To this end, for the sake of the argument, let us reconsider the limiting process from the preceding section by assuming that the $f_k$'s are actual renormalized solutions of the Vlasov-Maxwell-Boltzmann system (without defect). Of course, appropriate truncations of the system would have to be used here to be fully rigorous. In this case, the non-negative defect measure $\nu_\beta$ is obtained exactly as the defect in the limit \eqref{limit defect}. That is to say, it holds that, up to extraction,
\begin{equation*}
	\lim_{k\rightarrow\infty}\int_{[0,\infty)\times\mathbb{R}^3\times\mathbb{R}^3}
	\left(\beta'(f_k)Q(f_k,f_k)-\beta'(f)Q(f,f)\right)\varphi dtdxdv
	=
	\int_{[0,\infty)\times\mathbb{R}^3\times\mathbb{R}^3}
	\varphi \nu_\beta(dtdxdv),
\end{equation*}
for every $\varphi\in C_0^\infty\left([0,\infty)\times\mathbb{R}^3\times\mathbb{R}^3\right)$. There is another non-negative Radon measure $\lambda$ on $[0,\infty)\times\mathbb{R}^3$, which is of interest and arises as the defect in the entropy dissipation. It is thus simply defined, up to extraction, by
\begin{equation*}
	\lim_{k\rightarrow\infty}\int_{[0,\infty)\times\mathbb{R}^3}
	\left(D(f_k)-D(f)\right)\varphi dtdx
	=
	\int_{[0,\infty)\times\mathbb{R}^3}
	\varphi \lambda(dtdx),
\end{equation*}
for every $\varphi\in C_0^\infty\left([0,\infty)\times\mathbb{R}^3\right)$.

The following theorem relates $\nu_\beta$ to $\lambda$ and, thus, shows that $\nu_\beta$ can be controlled by the entropy dissipation only. The precise assumptions needed here on the cross-section are slightly more restrictive that those found in the preceding theorem and so, for the sake of simplicity, we only deal with inverse power potentials below.

\begin{theorem}
	Let $b(z,\sigma)$ be a collision kernel which derives from an inverse power potential. Then, the renormalized solutions of the Vlasov-Maxwell-Boltzmann system with a defect measure constructed in the preceding section can be shown to satisfy the following entropy inequality with a defect measure:
	\begin{equation*}
		\begin{aligned}
			\int_{\mathbb{R}^3\times\mathbb{R}^3} \left(f \log {f \over M} - f +M\right) (t) dxdv
			& + \frac 1{2} \int_{\mathbb{R}^3} \left( |E|^2+ |B|^2 \right)(t) dx
			+ \int_0^t\int_{\mathbb{R}^3}D(f)(s) dx ds + \lambda\left([0,t)\times\mathbb{R}^3\right)
			\\
			& \leq \int_{\mathbb{R}^3\times\mathbb{R}^3} \left( f^\mathrm{in}\log {f^\mathrm{in} \over M}
			- f^\mathrm{in}+M \right) dxdv
			+ \frac1{2}\int_{\mathbb{R}^3}\left(  |E^\mathrm{in}|^2+  |B^\mathrm{in}|^2\right)  dx,
		\end{aligned}
	\end{equation*}
	where, for any admissible renormalization $\beta(z)$, it holds that
	\begin{equation*}
		\int_{[0,\infty)\times\mathbb{R}^3\times\mathbb{R}^3}\varphi(t,x) d\nu_\beta(t,x,v)
		\leq
		2\left( \left\|\left[\beta'(z)+z\beta''(z)\right]^-\right\|_\infty + \left\|\beta'(z)\right\|_\infty\right)
		\int_{[0,\infty)\times\mathbb{R}^3} \varphi(t,x) d\lambda(t,x),
	\end{equation*}
	for every non-negative $\varphi\in C_0^\infty\left([0,\infty)\times\mathbb{R}^3\right)$.
\end{theorem}

The above characterization of $\nu_\beta$ by the entropy dissipation turns out to be crucial in the study of hydrodynamic limits of the Vlasov-Maxwell-Boltzmann system (see \cite{arsenio6}). Its proof is somewhat technical and we refer to \cite{arsenio6} for its full justification.




\bibliographystyle{plain}
\bibliography{biblio.bib}





\end{document}